\magnification=\magstep1
\baselineskip=15pt
\mathsurround=1pt

\font\Bbb=msbm10

\def\pf#1:{\noindent{\bf#1}:}
\def\arrow-under #1 {\buildrel #1 \over \longrightarrow} 
\def\Ñ{\kern1pt\vrule height2.4ptwidth3ptdepth-1.9pt\kern1pt}
\def\¾{\hfill \rlap{$\sqcup$}$\sqcap$\par\bigskip}

\def\Zü{\hbox{\Bbb Z}}
\def\Rå{\hbox{\Bbb R}}

\def\scrptC{{\cal C}}

\def\a{\alpha}
\def\b{\beta}
\def\g{\gamma}
\def\e{\varepsilon}
\def\h{\eta}
\def\p{\pi}
\def\f{\varphi}
\def\y{\psi}

\def\c{\raise 1pt\hbox{$\chi$}}
\def\D{\Delta}
\def\S{\Sigma}

\def\bdyÆ{\partial}
\def\infÔ{\infty}
\def\empÑ{\hbox{\Bbb ?}}

\def\homot¦{\simeq}
\def\iso¤{\approx}
\def\­{\neq}
\def\Î{\in}
\def\É{\cdots}
\def\²{\leq}
\def\³{\geq}
\def\:{\,\colon}
\def\{\hbox{$\lower1.5pt\hbox{$^1$}\mskip-3.5mu/\mskip-2mu_2$}}
\def\unionÚ{\cup}
\def\UnionÚ{\bigcup}
\def\intÛ{\cap}
\def\IntÛ{\bigcap}
\def\cross´{\times}
\def\lbr›{\{}
\def\rbrÍ{\}}
\def\res|{\rlap{\raise1pt\hbox{$|$}}\lower1pt\hbox{$|$}}
\def\to¬{\rightarrow}
\def\To¬{\longrightarrow}
\def\mpstoØ{\mapsto}
\def\incl«{\hookrightarrow}
\def\subsetofÞ{\subset}
\def\tild'  {\widetilde}

\def\dfn{\it}
\def\rel{\ {\rm rel\ }}
\def\Diff{Dif\!f}

\centerline{\bf Spaces of Incompressible Surfaces}
\medskip 
\centerline{Allen Hatcher}
\bigskip

The homotopy type of the space of PL homeomorphisms of a Haken 3\Ñmanifold was computed in [H1], and with the subsequent proof of the Smale conjecture in [H2], the computation carried over to diffeomorphisms as well. These results were also obtained independently by Ivanov [I1,I2]. The main step in the calculation in [H1], though not explicitly stated in these terms, was to show that the space of embeddings of an incompressible surface in a Haken 3\Ñmanifold has components which are contractible, except in a few special situations where the components have a very simple noncontractible homotopy type. The purpose of the present note is to give precise statements of these embedding results along with simplified proofs using ideas from [H3]. We also rederive the calculation of the homotopy type of the diffeomorphism group of a Haken 3\Ñmanifold.

\medskip 
Let $ M $ be an orientable compact connected irreducible 3-manifold, and let $ S $ be an incompressible surface in $ M $, by which we mean:
\smallskip
\item{---}$ S $ is a compact connected surface, not $ S^2 $, embedded in $ M $
properly, i.e., $ S \intÛ \bdyÆ M = \bdyÆ S $.
\item{---}The normal bundle of $ S $ in $ M $ is trivial. Since $ M $ is orientable, this
just means $ S $ is orientable.
\item{---}The inclusion $ S \incl« M $ induces an injective homomorphism on $ \p_1 $. 
\medskip\noindent 
We denote by $ E(S,M \rel \bdyÆ S) $ the space of smooth embeddings $ S \to¬ M $ agreeing with the given inclusion $ S \incl« M $ on $ \bdyÆ S $. The group $ \Diff(S \rel \bdyÆ S) $ of diffeomorphisms $ S \to¬ S $ restricting to the identity on $ \bdyÆ S $ acts freely on $ E(S,M \rel \bdyÆ S) $ by composition, with orbit space the space $ P(S,M \rel \bdyÆ S) $ of subsurfaces of $ M $ diffeomorphic to $ S $ by a diffeomorphism restricting to the identity on $ \bdyÆ S $. We think of points of $ P(S,M \rel \bdyÆ S) $ as ``positions" or ``placements" of $ S $ in $ M $. There is a fibration
$$
\Diff(S \rel \bdyÆ S) \To¬ E(S,M \rel \bdyÆ S) \To¬ P(S,M \rel \bdyÆ S)
$$
giving rise to a long exact sequence of homotopy groups. It is known that $ \p_i\Diff(S \rel \bdyÆ S) $ is zero for $ i > 0 $, except when $ S $ is the torus $ T^2 $ and $ i = 1 $, when the inclusion $ T^2 \incl« \Diff(T^2) $ as rotations induces an isomorphism on $ \p_1 $. So the higher homotopy groups of $ E(S,M \rel \bdyÆ S) $ and $ P(S,M \rel \bdyÆ S) $ are virtually identical. 

\filbreak
\proclaim Theorem 1. Let $ S $ be an incompressible surface in $ M $. Then:
\item{\rm (a)}$ \p_iP(S,M \rel \bdyÆ S) = 0 $ for all $ i > 0 $ unless $ \bdyÆ S = \empÑ $ and $ S $ is the fiber of a surface bundle structure on $ M $. In this exceptional case the inclusion $ S^1 \incl« P(S,M) $ as the fiber surfaces induces an isomorphism on $ \p_i $ for all $ i > 0 $. 
\item{\rm (b)}$ \p_iE(S,M \rel \bdyÆ S) = 0 $ for all $ i > 0 $ unless $ \bdyÆ S = \empÑ $ and $ S $ is either a torus or the fiber of a surface bundle structure on $ M $. In these exceptional cases $ \p_iE(S,M) = 0 $ for all $ i > 1 $. In the surface bundle case the inclusion of the subspace consisting of embeddings with image a fiber induces an isomorphism on $ \p_1 $. When $ S $ is a torus but not the fiber of a surface bundle structure, the inclusion of the subspace consisting of embeddings with image equal to the given $ S $ induces an isomorphism on $ \p_1 $.

Here it is understood that $ \p_iE(S,M \rel \bdyÆ S) $ and $ \p_iP(S,M \rel \bdyÆ S) $ are to be computed at the basepoint which is the given inclusion $ S \incl« M $. 

It is not hard to describe precisely what happens in the exceptional cases of (b). In the surface bundle case consider the exact sequence
$$
0 \To¬ \p_1\Diff(S) \To¬ \p_1E(S,M) \To¬ \p_1P(S,M) \arrow-under {\bdyÆ} \p_0\Diff(S)
$$
By part (a) we have $ \p_1P(S,M) \iso¤ \Zü $. The boundary map takes a generator of this $ \Zü $ to the monodromy diffeomorphism defining the surface bundle. If this monodromy has infinite order in $ \p_0\Diff(S) $ then the boundary map is injective so $ \p_1E(S,M) \iso¤ \p_1\Diff(S) $. If the monodromy has finite order in $ \p_0\Diff(S) $, it is isotopic to a periodic diffeomorphism of this order, as Nielsen showed. Hence $ M $ is Seifert-fibered with coherently oriented fibers, and there is an action of $ S^1 $ on $ M $ rotating circle fibers, and taking fibers of the surface bundle structure to fibers. The orbit of the given embedding of $ S $ under this action then generates a $ \Zü $ subgroup of $ \p_1E(S,M) $, which is all of $ \p_1E(S,M) $ unless $ S $ is a torus, in which case it is easy to see that $ \p_1E(S,M) \iso¤ \Zü \cross´ \Zü \cross´ \Zü $.

Using these results we can deduce:

\proclaim Theorem 2. If $ M $ is an orientable Haken manifold then $ \p_i\Diff(M \rel \bdyÆ M) = 0 $ for all $ i > 0 $ unless $ M $ is a closed Seifert manifold with coherently orientable fibers. In the latter case the inclusion $ S^1 \incl« \Diff(M) $ as rotations of the fibers induces an isomorphism on $ \p_i $ for all $ i > 0 $, except when $ M $ is the 3\Ñtorus, in which case the circle of rotations $ S^1 \incl« \Diff(M) $ is replaced by the 3-torus of rotations $ M \incl« \Diff(M) $.

\pf Proof: Suppose first that $ \bdyÆ M \­ \empÑ $, hence $ M $ is automatically Haken if it is irreducible. By the theory of Haken manifolds, there exists an incompressible surface $ S \subsetofÞ M $ with $ \bdyÆ S \­ \empÑ $, in fact with $ \bdyÆ S $ representing a nonzero class in $ H_1(\bdyÆ M) $. Consider the fibration
$$
\Diff(M \rel \bdyÆ M \unionÚ S) \To¬ \Diff(M \rel \bdyÆ M) \To¬ E(S,M \rel \bdyÆ S)
$$
The fiber can be identified with $ \Diff(M' \rel \bdyÆ M') $ where $ M' $ is the compact manifold, possibly disconnected, obtained by splitting $ M $ along $ S $. If we know the theorem holds for $ M' $, then from the long exact sequence of homotopy groups of the fibration, together with part (b) of the preceding theorem, we deduce that the theorem holds for $ M $. Again by Haken manifold theory, there is a finite sequence of such splitting operations reducing $ M $ to a disjoint union of balls. By the Smale conjecture the theorem holds for balls, so by induction the theorem holds for $ M $.

When $ M $ is a closed Haken manifold we again consider the fibration displayed above,
with
$ S
$ now a closed incompressible surface. If we are not in the exceptional cases that $
\p_1E(S,M)
$ is nonzero, described in the paragraph before Theorem 2, the arguments in the
preceding paragraph apply since $ M' $ is a nonclosed Haken manifold, for which the
theorem has already been proved. 

There remain the cases that $ \p_1E(S,M) \­ 0 $. 

\noindent({\it i}) If $ S $ is the fiber of a surface bundle structure on $ M $ but not a torus, the result follows from the remarks preceding Theorem 2, describing how a generator of $ \p_1E(S,M) \iso¤ \Zü $ is represented by the orbit of the inclusion $ S \incl« M $ under the $ S^1 $ action rotating fibers of the Seifert fibering.

\noindent({\it ii}) If $ S $ is a torus but not the fiber of a surface bundle, we have $ \p_1E(S,M) \iso¤ \Zü \cross´ \Zü $, the rotations of $ S $. Under the boundary map $ \p_1E(S,M) \to¬ \p_0\Diff(M \rel S) $ each loop of rotations goes to a Dehn twist on one side of $ S $ and its inverse twist on the other side. The behavior of such twists in the mapping class group of the manifold $ M' $ is well known; see e.g. [HM]. In particular, the boundary map is injective unless $ M $ is orientably Seifert-fibered with $ S $ a union of fibers, in which case the kernel of the boundary map is $ \Zü $ represented by rotations of the circle fibers.

\noindent({\it iii})  If $ S $ is the torus fiber of a surface bundle we have to
consider the phenomena in both ({\it i}) and ({\it ii}).  The monodromy diffeomorphism of
the torus fiber defining the bundle is an element of $ SL_2(\Zü) $. Note that a loop of
rotations of $ S $ as in ({\it ii}) lies in the kernel of the boundary map iff the
rotations are in the direction of a curve in $ S $ fixed by the monodromy. If the
monodromy is trivial, $ M $ is the 3\Ñtorus and we have contributions to $ \p_1\Diff(M)
$ from both ({\it i}) and ({\it ii}), so we get $ M \subsetofÞ \Diff(M) $ inducing an
isomorphism on $ \p_i $ for $ i > 0 $. If the monodromy is nontrivial and of finite
order, it has no real eigenvalues, so the boundary map in ({\it ii}) is injective and we
get only the $ S^1 \subsetofÞ \Diff(M) $ as in ({\it i}). If the monodromy is of
infinite order and has 1 as an eigenvalue, it is a Dehn twist of the torus fiber, so we
get $ \p_1\Diff(M) \iso¤ \Zü $ and this loop of diffeomorphisms is realized by rotating
fibers of a fibering of $ M $ by circles in the eigendirection in the torus fibers of $
M $. \¾

\medskip\noindent 
{\bf Proof of Theorem 1.}

We will first prove statement (b) when $ \bdyÆ S \­ \empÑ $, which suffices to deduce Theorem 2 in the nonclosed case. Then we will prove (a), and finally the remaining cases of (b).

Let $ f_t \: S \to¬ M $, $ t \Î D^i $, be a family of embeddings representing an element of $ \p_iE(S,M \rel \bdyÆ S) $. We assume $ i > 0 $, so all the surfaces $ f_t(S) $ are isotopic $ \rel \bdyÆ S $ to the given inclusion $ S \incl« M $, hence are incompressible also. Let $ S \cross´ I $ be a collar on one side of $ S = S \cross´ \{ 0 \} $ in $ M $. By Sard's theorem, $ f_t(S) $ is transverse to $ S \cross´ \{ x \} $ for almost all $ x \Î I $. Transversality is preserved under small perturbations, so, since $ D^i $ is compact, we can choose a finite cover of $ D^i $ by open sets $ U_j $ such that $ f_t(S) $ is transverse to a slice $ S_j = S \cross´ \{ x_j \} \subsetofÞ S \cross´ I $ for all $ t \Î U_j $. Then $ f_t(S) \intÛ S_j $ is a finite collection of disjoint circles which vary by isotopy as $ t $ ranges over $ U_j $. The main work will be to deform the family $ f_t $ to eliminate these circles, for all $ t $ and $ j $ simultaneously (after choosing the slices $ S_j $ and the open sets $ U_j $ a little more carefully).
\medskip
\noindent{\it Step 1.} The aim here is to eliminate all the nullhomotopic circles. Let $ \scrptC^j_t $ be the collection of circles of $ f_t(S) \intÛ S_j $ which are homotopically trivial in $ M $ and hence bound disks in both $ f_t(S) $ and $ S_j $, by incompressibility. Let $ \scrptC_t = \UnionÚ_j\scrptC^j_t $, the union over those $ j $'s such that $ t \Î U_j $. We would like to construct a family of functions $ \f_t $ assigning a value $ \f_t(C) \Î (0,1) $ to each circle $ C \Î \scrptC_t $ such that:
\smallskip
\item{(1)}$ \f_t(C) $ varies continuously with $ t \Î U_j $ for each $ C \Î \scrptC^j_t $. 
\item{(2)}$ \f_t(C) < \f_t(C') $ whenever the disk in $ f_t(S) $ bounded by $ C $ is contained in the disk bounded by $ C' $. 
\item{(3)}$ \f_t $ is injective for each $ t $, so $ \f_t(C) \­ \f_t(C') $ if $ C $ and $ C' $ are distinct circles of $ \scrptC_t $. 
\smallskip\noindent
Achieving (1) and (2) is not hard. For example, one can take $ \f_t(C) $ to be the area of the disk in $ f_t(S) $ bounded by $ C $, with respect to some metric on $ M $. To achieve (3) takes more work. First replace each $ S_j $ by $ 2i +1 $ nearby slices $ S_{jk} = S \cross´ \{ x_{jk} \} $, so that each circle of $ f_t(S) \intÛ S_j $ is replaced by $ 2i +1 $ nearby circles of $ f_t(S) \intÛ S_{jk} $. Define $ \f_t $ on each of these new circles of $ f_t(S) \intÛ S_{jk} $ to have value near the value of the original $ \f_t $ on the nearby circle of $ f_t \intÛ S_j $, so that (1) and (2) are still satisfied for all the new circles. Perturb the new functions $ t \mpstoØ \f_t(C) $ so that the solution set of each equation $ \f_t(C) = \f_t(C') $ is a codimension-one submanifold of $ D^i $, and so that these codimension-one submanifolds have general position intersections with each other, and in particular so that at most $ i $ such equations are satisfied for each $ t $. Then for each $ t $ the perturbed $ \f_t $ is injective on the complement of a set of at most $ 2i $ circles. This means that if we delete from $ \scrptC_t $ those circles on which $ \f_t $ is not injective, there exists for each $ t $ a slice $ S_{jk} $ from which no circles have been deleted. This slice has the same property for nearby $ t $. Let $ U_{jk} $ be the subspace of $ U_j $ consisting of points $ t $ for which no circles of $ \scrptC_t $ in $ S_{jk} $ are deleted. Then the open cover $ \{ U_{jk} \} $ of $ D^i $ and the associated slices $ S_{jk} $ satisfy (1)-(3). We relabel these as $ U_j $ and $ S_j $.

We would like to deform the family $ f_t $, $ t \Î D^i $, so as to eliminate all the circles of $ \scrptC_t $ without introducing new circles of $ f_t(S) \intÛ S_j $ for $ t \Î U_j $. Consider first a fixed value of $ t $ and a circle $ C \Î \scrptC_t $ for which $ \f_t(C) $ is minimal. Thus $ C $ bounds a disk $ D \subsetofÞ f_t(S) $ disjoint from all other circles of $ \scrptC_t $. In particular, $ D \intÛ S_j = C $, where $ C \subsetofÞ S_j $. By incompressibility of $ S_j $, $ C $ also bounds a disk $ D_j \subsetofÞ S_j $, and the sphere $ D \unionÚ D_j $ bounds a ball $ B \subsetofÞ M $ since $ M $ is irreducible. We can isotope $ f_t $ to eliminate $ C $ from $ f_t(S) \intÛ S_j $ by isotoping $ D $ across $ B $ to $ D_j $, and slightly beyond. This isotopy extends to an ambient isotopy of $ M $ supported near $ B $ which also eliminates any circles of $ f_t(S) \intÛ S_j $ which happen to lie inside $ D_j $; we call such circles {\dfn secondary\/}, in contrast to $ C $ itself which we call {\dfn primary\/}. Since $ D $ and $ D_j $ are disjoint from $ S_k $ if $ k \­ j $ and $ t \Î U_k $, so is $ B $, so this isotopy leaves circles of $ f_t(S) \intÛ S_k $ unchanged if $ k \­ j $ and $ t \Î U_k $. Hence we can iterate the process, eliminating in turn each remaining circle of $ \scrptC_t $ with smallest $ \f_t $ value. Thus we construct an isotopy $ f_{tu} $, $ 0 \² u \² 1 $, of $ f_t = f_{t0} $ eliminating each primary circle $ C \Î \scrptC_t $ during the $ u $-interval $ [\f_t(C),\f_t(C) + \e] $, for some fixed $ \e $, along with any secondary circles associated to $ C $. 

This isotopy $ f_{tu} $ will not depend continuously on $ t $ since as $ t $ moves from $ U_j $ to the complement of $ U_j $, we suddenly stop performing the isotopies eliminating the circles of $ f_t(S) \intÛ S_j $. This problem is easy to correct by the following truncation process. For each $ U_j $ choose a map $ \y_j \: D^i \to¬ [0,1] $ which is $ 0 $ outside $ U_j $ and $ 1 $ inside a slightly smaller open set $ U'_j $ in $ U_j $, such that the $ U'_j $'s still cover $ D^i $. Then modify the construction of $ f_{tu} $ by performing the isotopies eliminating primary circles of $ f_t(S) \intÛ S_j $ only for $ u \² \y_j(t) $. As observed earlier, the isotopy eliminating a primary circle of $ f_t(S) \intÛ S_j $ does not affect circles of $ f_t(S) \intÛ S_k $ for $ k \­ j $ with $ t \Î U_k $, so truncating such an isotopy at $ u = \y_j(t) $ creates no problems for continuing the construction of $ f_{tu} $ for $ u > \y_j(t) $ to eliminate circles of $ f_t(S) \intÛ S_k $. 

There is one other reason why $ f_{tu} $, as described so far, may  not depend
continuously on $ t $, namely, there is a choice in the isotopy eliminating a primary
circle, and these choices need to be made continuously in $ t $. We may specify an
isotopy eliminating a circle $ C $ by the following process. First enlarge the disk $ D
$ slightly to a disk $ D' \subsetofÞ f_t(S) $, with $ \bdyÆ D' $ contained in a nearby
parallel copy $ S'_j $ of $ S_j $, then choose a collar on $ D' $ containing $ B $,
i.e., an embedding $ \c \: D' \cross´ I \to¬ M $ with $ \c\res| D' \cross´ \{ 0 \} $ the
identity and $ B \subsetofÞ \c(D' \cross´ [0,\]) \subsetofÞ \c(D' \cross´ I) \subsetofÞ
B' $ where $ B' $ is a ball constructed like $ B $ using $ S'_j $ in place of $ S_j $.
We may also assume $ \c(\bdyÆ D' \cross´ I) \subsetofÞ S'_j $. Then an isotopy
eliminating $ C $ is obtained by pushing points $ x \Î D' $ along the lines $ \c(\{ x \}
\cross´ I) $, with this motion damped down near $ \bdyÆ D' $. 

The space of such collars is contractible. Namely, use one collar $ \h $ to produce an isotopy $ h_s \: B' \to¬ B' $ moving $ \h(D' \cross´ \{ 0 \}) $ to $ \h(D' \cross´ \{ 1 \}) $ in the obvious way. Then for an arbitrary collar $ \c $, construct a deformation $ \c_s $ of $ \c $ consisting of $ h_s\c $ together with a portion of $ \h $. (To make $ \c_s $ smooth at the junction of these two pieces one can require all $ \c $'s to have the same derivative along $ D' \cross´ \{ 0 \} $.) Then $ \c_1 $ contains the collar $ \h $, so we can deform $ \c_1 $ to $ \h $ by gradually truncating the part outside $ \h $. 

To make the dependence of $ f_{tu} $ on $ u $ explicit we may proceed as follows. We may assume the truncation functions $ \y_j(t) $ are piecewise linear, and we may perturb the functions $ \f_t(C) $ to be piecewise linear functions of $ t $ also. Then we may choose a triangulation of $ D^i $ so that each simplex is contained in some $ U_j $ and so that the solution set of each equation $ \y_j(t) = \f_t(C) $ is a subcomplex of the triangulation. Now we construct the isotopies $ f_{tu} $ inductively over skeleta of the triangulation. Assume that $ f_{tu} $ has already been constructed over the boundary of a $ k $-simplex $ \D^k $, and assume that over $ \D^k $ itself we have already constructed $ f_{tu} $ for $ t \² \f_t(C) $, for some primary circle $ C \Î \scrptC^j_t $. If $ \y_j(t) \² \f_t(C) $ over $ \D^k $, the induction step is vacuous, so we may suppose $ \y_j(t) \³ \f_t(C) $ for $ t \Î \D^k $. By induction, for $ t \Î \bdyÆ \D^k $ we have already chosen collars $ \c_t $ for the disk $ D' $, as above, depending continuously on $ t $, and since the space of collars is contractible we can extend these collars over $ \D^k $. This allows the induction step to be completed, eliminating $ C $ over $ \D^k $, with the elimination isotopy truncated if $ \y_t $ so dictates.
\medskip 
\noindent{\it Step 2.} We show how to finish the proof of (b) if $ S $ is a disk. After Step 1, all the circles of $ f_t(S) \intÛ S_j $ have been eliminated for all $ t \Î U_j $ and all $ j $. By averaging the slices $ S_j = S \cross´ \{ x_j \} $ via a partition of unity subordinate to the $ U_j $'s we can choose a continuously varying slice $ S_t = S \cross´ \{ x_t \} $ disjoint from $ f_t(S)$; here we use the fact that if $ f_t(S) $ is disjoint from $ S \cross´ \{ x \} $ and from $ S \cross´ \{ y \} $ then it is disjoint from $ S \cross´ [x,y] $ since $ f_t(S) $ is connected and its boundary, which is non-empty, lies in $ S = S \cross´ \{ 0 \} $. Having $ f_t(S) $ disjoint from $ S_t $ for all $ t $, we can then by isotopy extension isotope the family $ f_t $ so that $ f_t(S) $ is disjoint from $ S \cross´ \{ 1 \} $ for all $ t $. 

The proof can now be completed as follows. The space $ E $ of embeddings $ S \cross´ I \to¬ M $ agreeing with the given embedding on $ S \cross´ \{ 1 \} \unionÚ \bdyÆ S \cross´ I $ fits into a fibration:
$$
\Diff(S \cross´ I \rel \bdyÆ ) \To¬ E \To¬ E(S, M - S \cross´ \{ 1 \} \rel \bdyÆ S)
$$
It is elementary that $ E $ has trivial homotopy groups since one can canonically isotope an embedding in $ E $ so that it equals the given embedding on a neighborhood of $ S \cross´ \{ 1 \} \unionÚ \bdyÆ S \cross´ I $, then gradually excise from the embedding everything but this neighborhood. The fiber $ \Diff(S \cross´ I \rel \bdyÆ ) $  also has trivial homotopy groups by the Smale conjecture since $ S $ is a disk. Thus $ \p_iE(S, M - S \cross´ \{ 1 \} \rel \bdyÆ S) $ vanishes for $ i > 0 $, and the first part of the proof shows that this group maps onto $ \p_iE(S, M \rel \bdyÆ S) $, so the latter group also vanishes. When $ S $ is a disk this argument also applies for $ i = 0 $ since $ \p_0E(S, M - S \cross´ \{ 1 \} \rel \bdyÆ S) = 0 $ by the irreducibility of $ M $. 
\medskip 
\noindent{\it Step 3} is to eliminate the remaining circles of $ f_t(S) \intÛ S_j $ for all $ t \Î U_j $ and all $ j $, in the case $ \bdyÆ S \­ \empÑ $. The role of the ball $ B $ in Step 1 will be played by what we may call a {\dfn pinched product\/}. This is obtained from a product $ W \cross´ I $, where $ W $ is a nonclosed compact orientable surface, by collapsing each segment $ \{ w \} \cross´ I $, $ w \Î \bdyÆ W $, to a point. Thus a pinched product $ P $ is a handlebody with its boundary decomposed as the union of two copies $ \bdyÆ_+P $ and $ \bdyÆ_-P $ of the surface $ W $, with $ \bdyÆ_+P $ and $ \bdyÆ_-P $ intersecting only in their common boundary, a ``corner" of $ \bdyÆ P $. 

If we can find a pinched product $ P \subsetofÞ M $ such that $ \bdyÆ_+P \subsetofÞ f_t(S) $ and $ \bdyÆ_-P \subsetofÞ S_j $, then by isotoping $ f_t $ by pushing $ \bdyÆ_+P $ across $ P $ to $ \bdyÆ_-P $, and slightly beyond, we eliminate the circles of $ \bdyÆ_+P \intÛ \bdyÆ_-P $ from $ f_t(S) \intÛ S_j $, as well as any other circles of $ f_t(S) \intÛ S_j $ which happen to lie in $ \bdyÆ_-P $. 

In order to locate such pinched products it is convenient to consider the covering space $ p \: (\tild'  M,\tild'  x_0) \to¬ (M,x_0) $ corresponding to the subgroup $ \p_1(S,x_0) $ of $ \p_1(M,x_0) $, with respect to a basepoint $ x_0 \Î S $. There is then a homeomorphic copy $ \tild' S $ of $ S $ in $ \tild'  M $ containing $ \tild'  x_0 $. We can associate to $ \tild'  M $ a graph $ T $ having a vertex for each component of $ \tild'  M - p^{-1}(S) $ and an edge for each component of $ p^{-1}(S) $. This graph $ T $ is in fact a tree. For suppose $ \g $ is a loop in $ T $ based at a point in the edge corresponding to $ \tild' S $. This can be lifted to a loop $ \tild' \g $ in $ \tild' M $ based at $ \tild' x_0 $. By the definition of $ \tild'  M $, the loop $ \tild' \g $ is homotopic to a loop in $ \tild' S $, and this homotopy projects to a homotopy from $ \g $ to a trivial loop.

In the same way, the surface $ S_j $ parallel to $ S $ determines a tree $ T_j $ canonically isomorphic to $ T $. The lift $ \tild' S $ lies in a component of $ \tild'  M - p^{-1}(S_j) $ corresponding to a base vertex of $ T_j $, and hence to a base vertex of $ T $ which is independent of $ j $. 

The family $ f_t(S) $ is an $ i $-parameter isotopy of $ S $, so there are lifts $ \tild'  f_t \: S \to¬ \tild'  M $ such that $ \tild'  f_t(S) $ forms an $ i $-parameter isotopy of $ \tild' S $. The lifts $ \tild'  f_t(S) $ and the parameter domain being compact, we can choose a vertex $ v $ of $ T $ farthest from the base vertex among all vertices corresponding to components of $ \tild' M - p^{-1}(S_j) $ which meet $ \tild'  f_t(S) $ as $ t $ varies over $ S^k $ and $ j $ varies arbitrarily. Let $ \tild' V_j $ be the closure of the component of $ \tild' M - p^{-1}(S_j) $ corresponding to $ v $, with $ \tild' S_j $ its boundary component in the direction of $ \tild' S $. Let $ \tild' \scrptC^j_t $ be the collection of components of $ \tild'  f_t(S) \intÛ \tild' V_j $ and let $ \scrptC^j_t $ be the (diffeomorphic) images of these components in $ f_t(S) $.  Let $ \scrptC_t = \UnionÚ_j\scrptC^j_t $, the union over $ j $ such that $ t \Î U_j $. 

\proclaim Lemma. For each surface $ C \Î \scrptC^j_t $ there is a pinched product $ P \subsetofÞ M $ with $ \bdyÆ_+P = C $ and $ \bdyÆ_-P \subsetofÞ S_j $. 

\pf Proof: This is proved in II.5 of [L], but let us sketch an argument which may be more direct. First observe that $ \p_1(\tild' V_j,\tild' S_j) = 0 $, otherwise $ \p_1(\tild' S_j) \to¬ \p_1(\tild' V_j) $ would not be surjective and hence $ \p_1(\tild' S) \to¬ \p_1(\tild'  M) $ could not be an isomorphism. (This uses the fact that the components of $ p^{-1}(S_j) $ are incompressible in $ \tild'  M $.) Next, let $ V_j $ be $ M $ split along $ S_j $, so we have a covering space $ p \: (\tild' V_j,\tild' S_j) \to¬ (V_j,S_j) $. Since the inclusion $ (C,\bdyÆ C) \incl« (V_j,S_j) $ lifts to $ (\tild' V_j,\tild' S_j) $, the fact that $ \p_1(\tild' V_j,\tild' S_j) = 0 $ implies that $ \p_1(C,\bdyÆ C) \to¬ \p_1(V_j,S_j) $ is zero. Thus there is a map $ F \: D^2 \to¬ V_j $ giving a homotopy from a path $ \a $ in $ C $ representing a nontrivial element of $ \p_1(C,\bdyÆ C) $ to a path $ \b $ in $ S_j $. We would like to improve $ F $ to be an embedding with $ F(D^2) \intÛ C = \a $. To do this, first perturb $ F $ to be transverse to $ C $. We can modify $ F $ to eliminate any circles of $ F^{-1}(C) $, using the fact that $ C $ is incompressible in $ V_j $. Arcs of $ F^{-1}(C) $, other than $ \a $, which are trivial in $ \p_1(C,\bdyÆ C) $ can be eliminated similarly. An outermost remaining arc of $ F^{-1}(C) $ can be used as a new $ \a $ for which $ F(D^2) \intÛ C = \a $. Now we apply the loop theorem to replace $ F $ by an embedding with the same properties. Having this embedded disk, we can isotope $ C $ by pushing $ \a $ across the disk, surgering $ C $ to a simpler surface $ C' $ which, by induction on the complexity of $ C $, splits off a pinched product $ P' $ from $ V_j $; the induction starts with the case that $ C $ is a disk, where incompressibility of $ S_j $ and irreducibility of $ M $ gives the result. We recover $ C $ from $ C' $ by adjoining a ``tunnel." If this tunnel lies outside $ P' $, then the tunnel enlarges $ P' $ to a pinched product $ P $ split off from $ V_j $ by $ C $, as desired. The other alternative, that the tunnel lies inside $ P' $, cannot occur since $ C $ would then be compressible. \¾

Continuing with the main line of the proof, we would like to choose functions $ \f_t \: \scrptC_t \to¬ (0,1) $ satisfying:

\item{(1)}$ \f_t(C) $ varies continuously with $ t \Î U_j $ for each $ C \Î \scrptC^j_t $. 
\item{(2)}$ \f_t(C) < \f_t(C') $ if $ C \subsetofÞ C' $. 
\item{(3)}$ \f_t $ is injective for each $ t $. 
\smallskip\noindent
As before, (1) and (2) are easy to arrange, and then (3) is achieved by replacing each $ S_j $ with a number of nearby copies of itself and rechoosing the cover $ \{ U_j \} $. 

Having functions $ \f_t $ satisfying (1)-(3), we follow the same scheme as in Step 1 to
construct an isotopy $ f_{tu} $ of $ f_t $ eliminating components $ C \Î \scrptC_t $
during the corresponding $ u $-intervals $ [\f_t(C),\f_t(C) + \e] $. Again there are
primary and secondary components, and only the primary components need be dealt with
explicitly.

 The end result of this family of isotopies $ f_{tu} $ is that the vertex $ v
$ is no longer among the vertices of $ T $ corresponding to components of $ \tild' M -
p^{-1}(S_j) $ meeting $ \tild'  f_t(S) $ for $ t \Î U_j $, and no new such vertices have
been introduced. So by iteration of the process we eventually reach the situation that $
f_t(S) $ is disjoint from $ S_j $ for all $ t \Î U_j $ and all $ j $. 
\medskip 
\noindent{\it Step 4.} We can now finish the proof of (b) in the case that $ \bdyÆ S \­ \empÑ $ by the argument in Step 2. The assumption that $ S $ was a disk rather than an arbitrary compact orientable surface with non-empty boundary was used in Step 2 only to deduce that $ \Diff(S \cross´ I \rel \bdyÆ ) $ had trivial homotopy groups, but the case $ S = D^2 $ suffices to show this in the present case since the handlebody $ S \cross´ I $ can be reduced to a ball by cutting along a collection of disjoint disks; see the proof of Theorem 2.
\medskip 
\noindent{\it Step 5.} Now we prove (a) when $ \bdyÆ S \­ \empÑ $. Steps 1 and 3 work equally well in this case, with images of embeddings instead of actual embeddings. The only modification needed in steps 2 and 4 is to show that $ P(S,M - S \cross´ \{ 1 \} \rel \bdyÆ S) $ has trivial $ \p_i $ for $ i > 0 $. The basepoint component of $ P(S,M - S \cross´ \{ 1 \} \rel \bdyÆ S) $ can be identified with the base space in the following fibration, where $ A = S \cross´ \{ 1 \} \unionÚ \bdyÆ S \cross´ I $.
$$
\Diff(S \cross´ I \rel A) \To¬ E(S \cross´ I,M \rel A) \To¬ E(S \cross´ I,M \rel A)/\Diff(S \cross´ I \rel A)
$$
The total space $ E(S \cross´ I,M \rel A) $ is contractible, as in Step 2. And the fiber has trivial homotopy groups, as one can see from the fibration
$$
\Diff(S \cross´ I \rel \bdyÆ(S \cross´ I)) \To¬ \Diff(S \cross´ I \rel A) \To¬ \Diff(S \cross´ \{ 0 \} \rel \bdyÆ S \cross´ \{ 0 \})
$$

\noindent{\it Step 6.} This is to prove (a) when $ \bdyÆ S = \empÑ $. Let $ \S_t $ be  a
family of surfaces representing an element of $ \p_iP(S,M) $. Steps 1 and 3 work  when $
\bdyÆ S = \empÑ $, so we may assume $ \S_t $ is disjoint from $ S_j $ for all $ t \Î U_j
$ and all $ j $. However, this does not imply $ \S_t $ is disjoint from the regions
between $ S_j $'s as it did when $ \bdyÆ S \­ \empÑ $. Instead, let us look at the lifts
$ \tild' \S_t $ to $ \tild' M $, which are well defined at least when the parameter
domain $ S^i $ is simply-connected, i.e., when $ i > 1 $. Let $ \tild' V_j $ and $
\tild' S_j $ be defined as in Step 3. Then $ \tild' \S_t \subsetofÞ \tild' V_j $ for $ t
\Î U_j $. If $ \tild' S $ is not contained in $ \tild' V_j $, then $ \tild' S $ and $
\tild' \S_t $ are disjoint for $ t \Î U_j $, hence, since they are isotopic, the region
between them in $ \tild' M $ is a product $ \tild' S \cross´ I $. The surface $
\tild' S_j $ lies in this region, hence must be compact too. Since $ \tild' S_j $ is an
incompressible surface in the interior of the product $ \tild' S \cross´ I $, the region
between $ \tild' S_j $ and $ \tild' \S_t $ must be a product, projecting
diffeomorphically to a product region between $ S_j $ and $ \S_t $. We can use such
products in place of the pinched products in Step 3, and thus isotope the family $ \S_t
$ to a new family for which $ \tild' V_j $ does contain $ \tild' S $. Then we can easily
isotope the family $ \S_t $ to be disjoint from $ S \cross´ \{ 1 \} $ for all $ t $ and
apply the argument in Step 5 to finish the proof that $ \p_iP(S,M) = 0 $ for $ i > 1 $.

When $ i = 1 $ we can replace the parameter domain $ S^1 $ by $ I $ and then the lifts $ \tild' \S_t $ exist, with $ \tild' \S_0 = \tild' S $. The difficulty is that $ \tild' \S_1 $ may be a different lift of $ S $ from $ \tild'  S $. If this happens, then in the notation of the preceding paragraph, $ \tild' V_j $ will not contain $ \tild' S $. In this case the product $ \tild' S \cross´ I $ will contain another lift of $ S $. The region between $ \tild' S $ and such a lift will also be a product $ \tild' S \cross´ I $, which means that $ M $ must be a surface bundle with $ S $ as fiber. Then we have a surjective homomorphism $ \p_1P(S,M) \to¬ \Zü $ which measures which lift of $ S $ to $ \tild'  M = S \cross´ \Rå $ the surface $ \tild' \S_1 $ is. On the kernel of this homomorphism the arguments of the preceding paragraph apply, and we deduce that this kernel is zero. Thus $ \p_1P(S,M) $ is $ \Zü $, represented by fibers of the surface bundle. When $ M $ is not a surface bundle with fiber $ S $ we must have $ \tild' \S_1 = \tild'  S $, and $ \p_1P(S,M) = 0 $. 

\noindent{\it Step 7.} When $ \bdyÆ S = \empÑ $ we can deduce (b) from (a) by looking at the long exact sequence of homotopy groups for the fibration $ \Diff(S) \to¬ E(S,M) \to¬ P(S,M) $. This is immediate except when $ M $ is a surface bundle with fiber $ S $. In this special case, exactness at $ \p_1E(S,M) $ implies that this fundamental group is represented by loops of embeddings $ S \incl« M $ with image a fiber. \¾

\bigskip\noindent
{\bf References}
\parindent=0pt
\medskip

[H1] A. Hatcher, Homeomorphisms of sufficiently large $P^2$\Ñirreducible 3\Ñmanifolds,
Topology 15 (1976), 343-347.

[H2] A. Hatcher, A proof of the Smale conjecture, Ann. of Math. 117 (1983), 553-607.

[H3] A. Hatcher, On the diffeomorphism group of $ S^1 \cross´ S^2 $, Proc.A.M.S. 83
(1981), 427-430.

[I1] N. Ivanov, Diffeomorphism groups of Waldhausen manifolds, J. Soviet Math. 12, No. 1 (1979), 115-118. Russian original: Research in Topology, II, Notes of LOMI scientific seminars, v. 66 (1976), 172-176.

[I2] N. Ivanov, Spaces of surfaces in Waldhausen manifolds, Preprint LOMI P-5-80 (1980), 
31 pp. (in Russian).

[L] F. Laudenbach, Topologie de la dimension trois: homotopie et isotopie, Ast\'erisque 12, Soc. Math. France, 1974.

[LZ] F. Laudenbach and O. Zambrano, Disjonction \`a parametres dans une vari\'et\'e de dimension 3, C.R.Acad.Sci.Paris Ser.A-B 288 (1979), A331-334.

\end